\pdfoutput=1
\RequirePackage{ifpdf}
\ifpdf 
\documentclass[pdftex]{sigma}
\else
\documentclass{sigma}
\fi

\newtheorem{conj}{\bf Conjecture}
\begin{document}

\allowdisplaybreaks

\renewcommand{\PaperNumber}{070}

\FirstPageHeading

\ShortArticleName{Ultradiscrete Painlev\'e VI with Parity Variables}

\ArticleName{Ultradiscrete Painlev\'e VI with Parity Variables}

\Author{Kouichi TAKEMURA and Terumitsu TSUTSUI}

\AuthorNameForHeading{K.~Takemura and T.~Tsutsui}

\Address{Department of Mathematics, Faculty of Science and Technology, Chuo University,\\
1-13-27 Kasuga, Bunkyo-ku Tokyo 112-8551, Japan}
\Email{\href{mailto:takemura@math.chuo-u.ac.jp}{takemura@math.chuo-u.ac.jp}}

\ArticleDates{Received July 15, 2013, in f\/inal form November 11, 2013; Published online November 19, 2013}

\Abstract{We introduce a~ultradiscretization with parity variables of the $q$-dif\/ference Painlev\'e VI
system of equations.
We show that ultradiscrete limit of Riccati-type solutions of $q$-Painlev\'e VI satisf\/ies the
ultradiscrete Painlev\'e VI system of equations with the parity variables, which is valid by using the
parity variables.
We study some solutions of the ultradiscrete Riccati-type equation and those of ultradiscrete Painlev\'e VI
equation.}

\Keywords{Painlev\'e equation; ultradiscrete; numerical solutions}

\Classification{39A13; 34M55; 37B15}

\section{Introduction} The Painlev\'e equations appear frequently in the problem of mathematical physics,
and they have extremely rich structures of mathematics~\cite{BB}.
The $q$-Painlev\'e equations are $q$-dif\/ference analogues of the Painlev\'e equations~\cite{RGH}, and
most of them have symmetry of af\/f\/ine Weyl groups, which play important roles to analyze integrable
systems~\cite{OGH,Sak0}.
On the other hand, cellular automaton has been studied actively and has been applied to vast areas of
science and technology.
Although some of cellular automaton describe complexity from their simple rule of evolution, some of them
have integrability~\cite{TTMS}.

The ultradiscrete Painlev\'e equations are systems of cellular automaton and they are obtained by suitable
limits (ultradiscretization) from the $q$-Painlev\'e equations~\cite{TTGOR}.
We now explain ultradiscretization.
We take positive variables $x$, $y$, $z$ and set $x=\exp (X/\varepsilon )$, $y=\exp (Y/\varepsilon )$ and
$z=\exp (Z/\varepsilon )$.
If the variables $x$, $y$, $z$ satisf\/ies $x+y=z$, then we have $\max (X,Y)=Z$ by the limit $\varepsilon \to
+0 $, which follows from the formula $ \lim\limits_{\varepsilon \to +0} \varepsilon \log (\exp
(X/\varepsilon )+\exp (Y/\varepsilon ))= \max (X,Y) $.
It is easy to conf\/irm that the relation $ xy=z$ (resp.\
 $ x/y=z$) corresponds to $X+Y=Z$ (resp.\
 $X-Y=Z$).
This procedure is sometimes called ultradiscretization.
The addition, the multiplication and the division correspond to taking the maximum, the addition and the
subtraction by the ultradiscretization.
However the subtraction $x-y$ is not well-behaved by the ultradiscretezation.
Moreover if some values of $x$, $y$, $z$ are negative, then the procedure does not work well.
To overcome these troubles, Satsuma and his collaborators~\cite{MIMS} introduced the ultradiscretization
with parity variables, and they obtained the ultradiscrete Painlev\'e II with parity variables~\cite{IS}.
Then the ultradiscrete Airy function with parity variables appears as a~special solution of the
ultradiscrete Painlev\'e II with parity variables~\cite{IS}.

The $q$-dif\/ference Painlev\'e VI system of equations was discovered by Jimbo and Sakai~\cite{JS}, and it
is written as
\begin{gather}
\frac{z(t)z(qt)}{b_3b_4}=\frac{(y(t)-t a_1)(y(t)-t a_2)}{(y(t)-a_3)(y(t)-a_4)},
\qquad
\frac{y(t)y(qt)}{a_3a_4}=\frac{(z(qt)-t b_1)(z(qt)-t b_2)}{(z(qt)-b_3)(z(qt)-b_4)},
\label{eq:yyzz}
\end{gather}
with the constraint $ b_1 b_2 a_3 a_4=q a_1 a_2 b_3 b_4 $.
The original Painlev\'e VI equation is recovered by the limit $q\to 1$ (see~\cite{JS}).

The $q$-Painlev\'e VI system has Riccati-type solutions in the special cases~\cite{JS}.
Namely, if $b_1 a_3=qa_1 b_3$, $b_2 a_4=a_2 b_4$ and the functions $y(t)$ and $z(t)$ satisfy the following
Riccati-type equation:
\begin{gather}
z(qt)=b_4\frac{y(t)-t a_2}{y(t)-a_4},
\qquad
y(qt)=a_3\frac{z(qt)-t b_1}{z(qt)-b_3},
\label{eq:yzqt}
\end{gather}
then the functions $y(t)$ and $z(t)$ satisfy the $q$-Painlev\'e VI system.
It is also known that the Riccati-type equation has solutions expressed by $q$-hypergeometric
functions~\cite{JS,Sak}.

In this article, we consider ultradiscretization of the $q$-Painlev\'e VI system with parity va\-riab\-les.
For each value $m$ of the independent variable, we associate the signs $\frak{y}_m, \frak{z}_m \in \{ \pm 1
\}$ and the amplitudes $Y_m, Z_m \in \mathbb {R} $.
Def\/ine a~parity function $S(\zeta )$ for a~sign variable $\zeta $ by $S(1) =0$ and $S(-1)=-\infty $.
We now introduce the ultradiscrete Painlev\'e VI system of equations with the variables $(\frak{y}_m,
Y_m)$ and $(\frak{z}_m,Z_m)$ by
\begin{gather}
\max\big[\max(A_1,A_2)+mQ+Y_{m}+B_3+B_4+S(\frak{y}_{m}),
\nonumber\\
\qquad
\max(2Y_{m},A_3+A_4)+Z_{m}+Z_{m+1}+S(\frak{z}_{m}\frak{z}_{m+1}),
\nonumber\\
\qquad
\max(A_3,A_4)+Y_{m}+Z_{m}+Z_{m+1}+S(-\frak{y}_{m}\frak{z}_{m}\frak{z}_{m+1})\big]
\nonumber\\
\quad=\max\big[\max(2mQ+A_1+A_2,2Y_{m})+B_3+B_4,
\nonumber\\
\qquad
\max(A_1,A_2)+mQ+Y_{m}+B_3+B_4+S(-\frak{y}_{m}),
\nonumber\\
\qquad
\max(2Y_{m},A_3+A_4)+Z_{m}+Z_{m+1}+S(-\frak{z}_{m}\frak{z}_{m+1}),
\nonumber\\
\qquad
\max(A_3,A_4)+Y_{m}+Z_{m}+Z_{m+1}+S(\frak{y}_{m}\frak{z}_{m}\frak{z}_{m+1})\big],
\label{eq:udP6zz}
\\
\max\big[\max(B_1,B_2)+m Q+Z_{m+1}+A_3+A_4+S(\frak{z}_{m+1}),
\nonumber
\\
\qquad
\max(2Z_{m+1},B_3+B_4)+Y_{m}+Y_{m+1}+S(\frak{y}_{m}\frak{y}_{m+1}),
\nonumber\\
\qquad
\max(B_3,B_4)+Y_{m}+Y_{m+1}+Z_{m+1}+S(-\frak{y}_{m}\frak{y}_{m+1}\frak{z}_{m+1})\big]
\nonumber\\
\quad
=\max\big[\max(2mQ+B_1+B_2,2Z_{m+1})+A_3+A_4,
\nonumber\\
\qquad
\max(B_1,B_2)+mQ+Z_{m+1}+A_3+A_4+S(-\frak{z}_{m+1}),
\nonumber\\
\qquad
\max(2Z_{m+1},B_3+B_4)+Y_{m}+Y_{m+1}+S(-\frak{y}_{m}\frak{y}_{m+1}),
\nonumber\\
\qquad
\max(B_3,B_4)+Y_{m}+Y_{m+1}+Z_{m+1}+S(\frak{y}_{m}\frak{y}_{m+1}\frak{z}_{m+1})\big],
\label{eq:udP6yy}
\end{gather}
with the constraint
\begin{gather}
B_1+B_2+A_3+A_4=Q+A_1+A_2+B_3+B_4.
\label{eq:udconstraint}
\end{gather}
On equations~\eqref{eq:udP6zz}, \eqref{eq:udP6yy} we ignore the terms containing $S(\zeta)=-\infty$ in the maximum.
Note that ultradiscretization of the $q$-Painlev\'e VI equation without parity variables was already
introduced by Ormerod~\cite{Or}, and it is recovered by f\/ixing the parity variables by $\frak{y}_m=\frak{z}_m=-1 $.
Equations~\eqref{eq:udP6zz}, \eqref{eq:udP6yy} are obtained by ultradiscretization of the $q$-Painlev\'e~VI system
under the condition $ a_i >0$ and $b_i>0$ $(i=1,2,3,4)$, where the condition is used in the process of
obtaining the original Painlev\'e~VI equation in~\cite{JS}, and we can also obtain the ultradiscrete
Painlev\'e~VI system of equations which admits the parities of the parameters (see
equations~\eqref{eq:udP6zzp}, \eqref{eq:udP6yyp}).
On the ultradiscrete Painlev\'e VI system of equations with parity variables, we have existence of the
solution of the initial value problem, although the uniqueness does not hold true.
We also ultradiscretize the Riccati-type equation with parity variables and show that any solutions of the
ultradiscrete Riccati-type equation satisfy the ultradiscrete Painlev\'e~VI system.
Here the ultradiscretization with parity variables is essential because the ultradiscrete Riccati-type
equation does not have any solutions in the case $\frak{y}_m=\frak{z}_m=-1 $.

We try to study the solutions of the ultradiscrete Riccati-type equation and those of the ultradiscrete
Painlev\'e VI equation.
We give examples of solutions which are described by piecewise-linear functions.
Based on numerical calculations, we present a~conjecture that the solutions are expressed as linear
functions if the independent variable $m$ is enough large (for details see Conjecture~\ref{conj:sol}).

This paper is organized as follows.
In Section~\ref{sec:UD}, we consider ultradiscretization of the Riccati-type equation and that of the
$q$-Painlev\'e VI equation.
We establish existence of the solution of the initial value problem for the ultradiscrete Painlev\'e VI
system of equations.
In Section~\ref{sec:UDRUDP6}, we show that any solutions of the ultradiscrete Riccati-type equation also
satisfy the ultradiscrete Painlev\'e VI equation.
In Section~\ref{sec:sol}, we investigate solutions of the ultradiscrete Riccati-type equation and those of
the ultradiscrete Painlev\'e VI equation.

\section{Ultradiscretization}\label{sec:UD}

\subsection{Ultradiscretization of the Riccati-type equation}

We consider ultradiscretization of the Riccati-type equation at f\/irst,
because the expression is simpler than those of the $q$-Painlev\'e VI equation.

To obtain the ultradiscrete limit, we set
\begin{gather}
t=q^m,
\qquad
q=e^{Q/\varepsilon},
\qquad
a_i=e^{A_i/\varepsilon}, \qquad b_i=e^{B_i/\varepsilon}, \qquad i=1,2,3,4,
\nonumber
\\
y\big(q^m\big)=(s(\frak{y}_{m})-s(-\frak{y}_{m}))e^{Y_m/\varepsilon},
\qquad
z(q^m)=(s(\frak{z}_{m})-s(-\frak{z}_{m}))e^{Z_m/\varepsilon},
\label{eq:limitsetting}
\end{gather}
and def\/ine the parity functions $s(\zeta )$ and $S(\zeta)$ by
\begin{gather*}
s(\zeta)=\begin{cases}
1,&\zeta=+1,
\\
0,&\zeta=-1,
\end{cases}
\qquad S(\zeta)=\begin{cases}
0,&\zeta=+1,
\\
-\infty,&\zeta=-1.
\end{cases}
\qquad
\end{gather*}
Then we have $s(\zeta )= e^{ S(\zeta)/\varepsilon }$ for $\varepsilon >0$.
We substitute equation~\eqref{eq:limitsetting} into the equation
\begin{gather*}
z(qt)(y(t)-a_4)=b_4(y(t)-t a_2),
\end{gather*}
which is equivalent to f\/irst equation of~\eqref{eq:yzqt}, and transpose the terms to disappear the minus
signs.
Then we have
\begin{gather*}
s(\frak{z}_{m+1})e^{Z_{m+1}/\varepsilon}s(\frak{y}_{m})e^{Y_{m}/\varepsilon}+s(-\frak{z}_{m+1})e^{Z_{m+1}
/\varepsilon}\big\{s(-\frak{y}_{m})e^{Y_{m}/\varepsilon}+e^{A_4/\varepsilon}\big\}
\\
\qquad
{}+e^{B_4/\varepsilon}\big\{s(-\frak{y}_{m})e^{Y_{m}/\varepsilon}+e^{(mQ+A_2)/\varepsilon}\big\}=e^{B_4/\varepsilon}
s(\frak{y}_{m})e^{Y_{m}/\varepsilon}
\\
\qquad
{}+s(\frak{z}_{m+1})e^{Z_{m+1}/\varepsilon}\big\{s(-\frak{y}_{m})e^{Y_{m}/\varepsilon}+e^{A_4/\varepsilon}\big\}
+s(-\frak{z}_{m+1})e^{Z_{m+1}/\varepsilon}s(\frak{y}_{m})e^{Y_{m}/\varepsilon}.
\end{gather*}
By using the formula $ s(\frak{y})s(\frak{z})+ s(-\frak{y})s(-\frak{z})= s(\frak{y} \frak{z})$
and taking the limit $\varepsilon \to +0$, we have
\begin{gather}
\max\big[mQ+A_2+B_4,Z_{m+1}+A_4+S(-\frak{z}_{m+1}),
\nonumber
\\
\qquad
Y_m+B_4+S(-\frak{y}_{m}),Y_m+Z_{m+1}+S(\frak{y}_{m}\frak{z}_{m+1})\big]
\nonumber
\\
\quad
=\max\big[Z_{m+1}+A_4+S(\frak{z}_{m+1}),Y_m+B_4+S(\frak{y}_{m}),Y_m+Z_{m+1}+S(-\frak{y}_{m}\frak{z}_{m+1})\big].
\label{eq:udRiccati2}
\end{gather}
It follows from the second equation of~\eqref{eq:yzqt} that
\begin{gather}
\max\big[mQ+A_3+B_1,Y_{m+1}+B_3+S(-\frak{y}_{m+1}),
\nonumber\\
\qquad
Z_{m+1}+A_3+S(-\frak{z}_{m+1}),Y_{m+1}+Z_{m+1}+S(\frak{y}_{m+1}\frak{z}_{m+1})\big]
\nonumber\\
\quad
=\max\big[Z_{m+1}+A_3+S(\frak{z}_{m+1}),Y_{m+1}+B_3+S(\frak{y}_{m+1}),
\nonumber\\
\qquad
Y_{m+1}+Z_{m+1}+S(-\frak{y}_{m+1}\frak{z}_{m+1})\big].
\label{eq:udRiccati1}
\end{gather}
We call equations~\eqref{eq:udRiccati2}, \eqref{eq:udRiccati1} the ultradiscrete Riccati-type equation with parity
variables.
By the ultradiscrete limit, the conditions $b_1 a_3=qa_1 b_3$, $b_2 a_4=a_2 b_4$ correspond to $B_1+A_3
=Q+A_1 +B_3$ and $B_2+A_4=A_2+B_4$.

We write the equations for the amplitude variables by f\/ixing the parity variables.
Equation~\eqref{eq:udRiccati2} for each case is written as
\begin{gather}
\frak{z}_{m+1}=1,\quad \frak{y}_{m}=1\qquad\Rightarrow
\nonumber
\\ \qquad\max(mQ+A_2+B_4,Y_m+Z_{m+1})=\max(Z_{m+1}+A_4,Y_m+B_4),\label{eq:udR2++}
\\
\frak{z}_{m+1}=1,\quad \frak{y}_{m}=-1\qquad\Rightarrow
\nonumber
\\ \qquad\max(mQ+A_2,Y_m)+B_4=Z_{m+1}+\max(Y_m,A_4),
\nonumber
\\
\frak{z}_{m+1}=-1,\quad \frak{y}_{m}=1\qquad\Rightarrow
\nonumber
\\
\qquad
\max(mQ+A_2+B_4,Z_{m+1}+A_4)=Y_m+\max(Z_{m+1},B_4),
\nonumber
\\
\frak{z}_{m+1}=-1,\quad \frak{y}_{m}=-1\qquad\Rightarrow
\nonumber
\\
\qquad
-\infty=\max(mQ+A_2+B_4,Y_m+Z_{m+1},Z_{m+1}+A_4,Y_m+B_4).
\nonumber
\end{gather}
There is no solution in the case $\frak{z}_{m+1}=\frak{y}_{m}=-1$.
Equation~\eqref{eq:udRiccati1} for each case is written as
\begin{gather}
\frak{z}_{m+1}=1,\quad \frak{y}_{m+1}=1\qquad\Rightarrow
\nonumber\\ \qquad\max(mQ+A_3+B_1,Y_{m+1}+Z_{m+1})=\max(Z_{m+1}+A_3,Y_{m+1}+B_3),
\label{eq:udR1++}
\\
\frak{z}_{m+1}=1,\quad \frak{y}_{m+1}=-1\qquad\Rightarrow
\nonumber\\ \qquad\max(mQ+A_3+B_1,Y_{m+1}+B_3)=Z_{m+1}+\max(Y_{m+1},A_3),
\nonumber\\
\frak{z}_{m+1}=-1,\quad \frak{y}_{m+1}=1\qquad\Rightarrow
\nonumber\\ \qquad\max(mQ+B_1,Z_{m+1})+A_3=Y_{m+1}+\max(Z_{m+1},B_3),
\nonumber\\
\frak{z}_{m+1}=-1,\quad \frak{y}_{m+1}=-1\qquad\Rightarrow
\nonumber\\
\qquad
-\infty=\max(mQ+A_3+B_1,Y_{m+1}+Z_{m+1},Z_{m+1}+A_3,Y_{m+1}+B_3).
\nonumber
\end{gather}
There is no solution in the case $\frak{z}_{m+1}=\frak{y}_{m+1}=-1$.

\subsection[Ultradiscretization of the $q$-Painlev\'e VI equation]{Ultradiscretization of the
$\boldsymbol{q}$-Painlev\'e VI equation} We can obtain similarly the ultradiscrete limit with parity
variables of the $q$-Painlev\'e VI equation.
We take a~limit of equation~\eqref{eq:yyzz} as $\varepsilon \to +0 $ by setting the values as
equation~\eqref{eq:limitsetting}.
By using the formulae
\begin{gather*}
s(\frak{y})+s(-\frak{y})=1,
\qquad
s(\frak{y})s(-\frak{y})=0,
\\
s(\frak{y})^2=s(\frak{y}),
\qquad
s(\frak{y})s(\frak{z})+s(-\frak{y})s(-\frak{z})=s(\frak{y}\frak{z}),
\\
s(\frak{y})s(\frak{z})s(\frak{w})+s(-\frak{y})s(-\frak{z})s(\frak{w})
+s(\frak{y})s(-\frak{z})s(-\frak{w})+s(-\frak{y})s(\frak{z})s(-\frak{w})=s(\frak{y}\frak{z}\frak{w}),
\end{gather*}
we obtain the ultradiscrete Painlev\'e VI equation with parity variables (i.e.\
equations~\eqref{eq:udP6zz}, \eqref{eq:udP6yy}).
In~\cite{Tsu}, another form of equations~\eqref{eq:udP6zz}, \eqref{eq:udP6yy} is derived with the details.
By the ultradiscrete limit, the constraint $ b_1 b_2 a_3 a_4=q a_1 a_2 b_3 b_4 $ of the $q$-Painlev\'e VI
system corresponds to equation~\eqref{eq:udconstraint}.

We write the equations for the amplitude variables by f\/ixing the parity variables.
Equation~\eqref{eq:udP6zz} for each case is written as
\begin{gather}
\frak{y}_{m}=1,\quad \frak{z}_{m}\frak{z}_{m+1}=1\qquad\Rightarrow
\nonumber
\\
\max(\max(2Y_{m},A_3+A_4)+Z_m+Z_{m+1},\max(A_1,A_2)+mQ+Y_{m}+B_3+B_4)
\nonumber
\\
=\max(\max(2mQ+A_1+A_2,2Y_{m})+B_3+B_4,\max(A_3,A_4)+Y_m+Z_m+Z_{m+1}),
\label{eq:om+zmzm1+}
\\
\frak{y}_{m}=1,\quad \frak{z}_{m}\frak{z}_{m+1}=-1\qquad\Rightarrow
\nonumber
\\
\max(\max(A_3,A_4)+Y_m+Z_m+Z_{m+1},\max(A_1,A_2)+mQ+Y_{m}+B_3+B_4)
\nonumber
\\
=\max(\max(2mQ+A_1+A_2,2Y_{m})+B_3+B_4,\max(2Y_{m},A_3+A_4)+Z_m+Z_{m+1}),
\label{eq:om+zmzm1-}
\\
\frak{y}_{m}=-1, \quad \frak{z}_{m}\frak{z}_{m+1}=1\qquad\Rightarrow
\nonumber\\
Z_m+Z_{m+1}+\max(A_3,Y_m)+\max(A_4,Y_m)
\nonumber\\
\qquad=B_3+B_4+\max(mQ+A_1,Y_m)+\max(mQ+A_2,Y_m),
\label{eq:om-zmzm1+}
\\
\frak{y}_{m}=-1,\quad \frak{z}_{m}\frak{z}_{m+1}=-1\qquad\Rightarrow
\nonumber\\
-\infty=\max(\max(2Y_{m},A_3+A_4)+Z_m+Z_{m+1},\max(A_1,A_2)+mQ+Y_{m}+B_3+B_4,
\nonumber
\\ \qquad\max(2mQ+A_1+A_2,2Y_{m})+B_3+B_4,\max(A_3,A_4)+Y_m+Z_m+Z_{m+1})).
\nonumber
\end{gather}
In the case $\frak{y}_{m}=-1$, $\frak{z}_{m} \frak{z}_{m+1}=1$, we used the associativity of the maximum
and the addition,~i.e.\
\begin{gather}
\max(X_1+W_1,X_2+W_1,X_1+W_2,X_2+W_2)=\max(X_1,X_2)+\max(W_1,W_2).
\label{eq:ass}
\end{gather}
There is no solution in the case $\frak{y}_{m}=\frak{z}_{m} \frak{z}_{m+1}=-1$.
On the other hand, equation~\eqref{eq:udP6yy} for each case is written as
\begin{gather}
\frak{z}_{m+1}=1,\quad \frak{y}_{m}\frak{y}_{m+1}=1\qquad\Rightarrow
\nonumber
\\
\max(\max(2Z_{m+1},B_3+B_4)+Y_m+Y_{m+1},\max(B_1,B_2)+m Q+Z_{m+1}+A_3+A_4)
\nonumber
\\
=\max(\max(2mQ+B_1+B_2,2Z_{m+1})+A_3+A_4,\max(B_3,B_4)\!+\!Y_m\!+\!Y_{m+1}\!+\!Z_{m+1}),
\label{eq:zm1+omom1+}
\\
\frak{z}_{m+1}=1, \quad \frak{y}_{m}\frak{y}_{m+1}=-1\qquad\Rightarrow
\nonumber\\
\max(\max(B_3,B_4)+Y_m+Y_{m+1}+Z_{m+1},\max(B_1,B_2)+m Q+A_3+A_4+Z_{m+1})
\nonumber\\
=\max(\max(2mQ+B_1+B_2,2Z_{m+1})+A_3+A_4,\max(2Z_{m+1},B_3+B_4)+Y_m+Y_{m+1}),
\nonumber\\
\frak{z}_{m+1}=-1,\quad \frak{y}_{m}\frak{y}_{m+1}=1\qquad\Rightarrow
\nonumber\\
Y_m+Y_{m+1}+\max(B_3,Z_{m+1})+\max(B_4,Z_{m+1})
\nonumber\\
=A_3+A_4+\max(mQ+B_1,Z_{m+1})+\max(mQ+B_2,Z_{m+1}),
\label{eq:om-ymym1+}
\\
\frak{z}_{m+1}=-1, \quad \frak{y}_{m}\frak{y}_{m+1}=-1\qquad\Rightarrow
\nonumber
\\
-\infty=\max(\max(2Z_{m+1},B_3+B_4)+Y_m+Y_{m+1},\max(B_1,B_2)+m Q+Z_{m+1}+A_3+A_4,
\nonumber
\\ \max(2mQ+B_1+B_2,2Z_{m+1})+A_3+A_4,\max(B_3,B_4)+Y_m+Y_{m+1}+Z_{m+1})).
\nonumber
\end{gather}
There is no solution in the case $\frak{z}_{m+1}=\frak{y}_{m} \frak{y}_{m+1}=-1$.

We can also obtain the ultradiscrete Painlev\'e VI system of equations which admits the parities of
parameters.
Set
\begin{gather*}
a_i=(s(\frak{a}_{i})-s(-\frak{a}_{i}))e^{A_i/\varepsilon}, \qquad b_i=(s(\frak{b}_{i})-s(-\frak{b}_{i}
))e^{B_i/\varepsilon}, \qquad i=1,2,3,4,
\end{gather*}
in addition to equation~\eqref{eq:limitsetting}, where $\frak{a}_{i}, \frak{b}_{i} \in \{\pm1\}$ represent
the signs of $a_i$, $b_i$ $(i=1,2,3,4)$.
We take the ultradiscrete limit ($\varepsilon \to +0$).
Then we have
\begin{gather}
\max\big[2mQ+A_1+A_2+B_3+B_4+S(-\frak{a}_{1}\frak{a}_{2}\frak{b}_{3}\frak{b}_{4}),
\nonumber\\ \qquad
2Y_{m}+B_3+B_4+S(-\frak{b}_{3}\frak{b}_{4}),
\nonumber\\ \qquad
Y_{m}+mQ+A_1+B_3+B_4+S(\frak{a}_{1}\frak{b}_{3}\frak{b}_{4}\frak{y}_{m}),
\nonumber\\ \qquad
Y_{m}+mQ+A_2+B_3+B_4+S(\frak{a}_{2}\frak{b}_{3}\frak{b}_{4}\frak{y}_{m}),
\nonumber\\ \qquad
2Y_{m}+Z_{m}+Z_{m+1}+S(\frak{z}_{m}\frak{z}_{m+1}),
\nonumber\\ \qquad
Z_{m}+Z_{m+1}+A_3+A_4+S(\frak{a}_{3}\frak{a}_{4}\frak{z}_{m}\frak{z}_{m+1}),
\nonumber\\ \qquad
Y_{m}+Z_{m}+Z_{m+1}+A_3+S(-\frak{a}_{3}\frak{y}_{m}\frak{z}_{m}\frak{z}_{m+1}),
\nonumber\\ \qquad
Y_{m}+Z_{m}+Z_{m+1}+A_4+S(-\frak{a}_{4}\frak{y}_{m}\frak{z}_{m}\frak{z}_{m+1})\big]
\nonumber\\
\quad
=\max\big[2mQ+A_1+A_2+B_3+B_4+S(\frak{a}_{1}\frak{a}_{2}\frak{b}_{3}\frak{b}_{4}),
\nonumber\\ \qquad
2Y_{m}+B_3+B_4+S(\frak{b}_{3}\frak{b}_{4}),
\nonumber\\ \qquad
Y_{m}+mQ+A_1+B_3+B_4+S(-\frak{a}_{1}\frak{b}_{3}\frak{b}_{4}\frak{y}_{m}),
\nonumber\\ \qquad
Y_{m}+mQ+A_2+B_3+B_4+S(-\frak{a}_{2}\frak{b}_{3}\frak{b}_{4}\frak{y}_{m}),
\nonumber\\ \qquad
2Y_{m}+Z_{m}+Z_{m+1}+S(-\frak{z}_{m}\frak{z}_{m+1}),
\nonumber\\ \qquad
Z_{m}+Z_{m+1}+A_3+A_4+S(-\frak{a}_{3}\frak{a}_{4}\frak{z}_{m}\frak{z}_{m+1}),
\nonumber\\ \qquad
Y_{m}+Z_{m}+Z_{m+1}+A_3+S(\frak{a}_{3}\frak{y}_{m}\frak{z}_{m}\frak{z}_{m+1}),
\nonumber\\ \qquad
Y_{m}+Z_{m}+Z_{m+1}+A_4+S(\frak{a}_{4}\frak{y}_{m}\frak{z}_{m}\frak{z}_{m+1})\big],\label{eq:udP6zzp}
\\
\max\big[2mQ+A_3+A_4+B_1+B_2+S(-\frak{a}_{3}\frak{a}_{4}\frak{b}_{1}\frak{b}_{2}),
\nonumber\\ \qquad
2Z_{m+1}+A_3+A_4+S(-\frak{a}_{3}\frak{a}_{4}),
\nonumber\\ \qquad
Z_{m+1}+m Q+A_3+A_4+B_1+S(\frak{a}_{3}\frak{a}_{4}\frak{b}_{1}\frak{z}_{m+1}),
\nonumber\\ \qquad
Z_{m+1}+m Q+A_3+A_4+B_2+S(\frak{a}_{3}\frak{a}_{4}\frak{b}_{2}\frak{z}_{m+1}),
\nonumber\\ \qquad
2Z_{m+1}+Y_{m}+Y_{m+1}+S(\frak{y}_{m}\frak{y}_{m+1}),
\nonumber\\ \qquad
Y_{m}+Y_{m+1}+B_3+B_4+S(\frak{b}_{3}\frak{b}_{4}\frak{y}_{m}\frak{y}_{m+1}),
\nonumber\\ \qquad
Y_{m}+Y_{m+1}+Z_{m+1}+B_3+S(-\frak{b}_{3}\frak{y}_{m}\frak{y}_{m+1}\frak{z}_{m+1}),
\nonumber\\ \qquad
Y_{m}+Y_{m+1}+Z_{m+1}+B_4+S(-\frak{b}_{4}\frak{y}_{m}\frak{y}_{m+1}\frak{z}_{m+1})\big]
\nonumber\\
\quad
=\max\big[2mQ+A_3+A_4+B_1+B_2+S(\frak{a}_{3}\frak{a}_{4}\frak{b}_{1}\frak{b}_{2}),
\nonumber\\ \qquad
2Z_{m+1}+A_3+A_4+S(\frak{a}_{3}\frak{a}_{4}),
\nonumber\\ \qquad
Z_{m+1}+m Q+A_3+A_4+B_1+S(-\frak{a}_{3}\frak{a}_{4}\frak{b}_{1}\frak{z}_{m+1}),
\nonumber\\ \qquad
Z_{m+1}+m Q+A_3+A_4+B_2+S(-\frak{a}_{3}\frak{a}_{4}\frak{b}_{2}\frak{z}_{m+1}),
\nonumber\\ \qquad
2Z_{m+1}+Y_{m}+Y_{m+1}+S(-\frak{y}_{m}\frak{y}_{m+1}),
\nonumber\\ \qquad
Y_{m}+Y_{m+1}+B_3+B_4+S(-\frak{b}_{3}\frak{b}_{4}\frak{y}_{m}\frak{y}_{m+1}),
\nonumber\\ \qquad
Y_{m}+Y_{m+1}+Z_{m+1}+B_3+S(\frak{b}_{3}\frak{y}_{m}\frak{y}_{m+1}\frak{z}_{m+1}),
\nonumber\\ \qquad
Y_{m}+Y_{m+1}+Z_{m+1}+B_4+S(\frak{b}_{4}\frak{y}_{m}\frak{y}_{m+1}\frak{z}_{m+1})\big],\label{eq:udP6yyp}
\end{gather}
with the constraint
\begin{gather*}
B_1+B_2+A_3+A_4=Q+A_1+A_2+B_3+B_4,
\qquad
\frak{a}_{1}\frak{a}_{2}\frak{a}_{3}\frak{a}_{4}=\frak{b}_{1}\frak{b}_{2}\frak{b}_{3}\frak{b}_{4}.
\end{gather*}
By setting $\frak{a}_{i} =\frak{b}_{i} =+1$ $(i=1,2,3,4)$, we recover
equations~\eqref{eq:udP6zz}, \eqref{eq:udP6yy}.
In the rest of the paper, we consider the case $\frak{a}_{i} =\frak{b}_{i} =+1$ $(i=1,2,3,4)$, i.e.\
equations~\eqref{eq:udP6zz}, \eqref{eq:udP6yy} for simplicity.

The ultradiscrete Painlev\'e VI equation with parity variables has solutions for any give initial values.
Namely we have the following proposition:
\begin{proposition}
Let $n_o \in \mathbb {Z} $, $\tilde{\frak{y} }_o, \tilde{\frak{z} }_o \in \{ \pm 1 \}$ and $y_o, z_o \in
\mathbb {R} $.
Then there exists a~solution $(\frak{y}_n, Y_n)$, $(\frak{z}_n, Z_n)$ $(n\in \mathbb {Z} )$ of
equations~\eqref{eq:udP6zz}, \eqref{eq:udP6yy} which satisfies the condition $(\frak{y}_{n_o}, Y_{n_o} )=
(\tilde{\frak{y} }_o, y_o)$ and $(\frak{z}_{n_o}, Z_{n_o} )= ( \tilde{\frak{z} }_o, z_o )$.
\end{proposition}
\begin{proof}
We show that, if the values $(\frak{z}_m, Z_m )$ and $(\frak{y}_m, Y_m )$ are f\/ixed, then there exists
$(\frak{z}_{m+1}$, $Z_{m+1})$ such that equation~\eqref{eq:udP6zz} is satisf\/ied.

If $\frak{y}_m =-1$, then we have $\frak{z}_{m+1}=\frak{z}_m $ and $Z_{m+1}$ is determined by
equation~\eqref{eq:om-zmzm1+}.

Assume that $\frak{y}_m =1$.
Write $U= \max( 2Y_{m},A_3+A_4)$, $U'=\max(A_3,A_4) +Y_m$, $V= \max( 2mQ+ A_1+A_2, 2Y_{m})$ and $V'=
\max(A_1,A_2) +mQ +Y_{m} $.
If $U\geq U'$ and $V\geq V'$ (resp.\
{}$U<U'$ and $V<V'$), then equation~\eqref{eq:udP6zz} is satisf\/ied by setting $\frak{z}_{m+1}=\frak{z}_m $
and $Z_{m+1}= V -U +B_3+B_4 -Z_m$ (resp.\
{}$Z_{m+1}= V' -U' +B_3+B_4 -Z_m$) (see equation~\eqref{eq:om+zmzm1+}).
If $U\geq U'$ and $V< V'$ (resp.\
{}$U<U'$ and $V\geq V'$), then equation~\eqref{eq:udP6zz} is satisf\/ied by setting $\frak{z}_{m+1}=-\frak{z}_m
$ and $Z_{m+1}= V' -U +B_3+B_4 -Z_m$ (resp.\
{}$Z_{m+1}= V -U' +B_3+B_4 -Z_m$) (see equation~\eqref{eq:om+zmzm1-}).
Similarly there exists $(\frak{y}_{m+1}, Y_{m+1} )$ such that equation~\eqref{eq:udP6yy} is satisf\/ied while
$(\frak{y}_m, Y_m )$ and $(\frak{z}_{m+1}, Z_{m+1} )$ are f\/ixed.

We proceed it for $m=n_o, n_o+1, \dots$, and a~solution $(\frak{y}_n, Y_n)$, $(\frak{z}_n, Z_n)$
$(n>n_o)$ is obtained.
By applying a~similar procedure to equation~\eqref{eq:udP6yy} for $n=n_o-1$, equation~\eqref{eq:udP6zz} for $n=n_o-1$
and so on, a~solution $(\frak{y}_n, Y_n)$, $(\frak{z}_n, Z_n)$ $(n<n_o )$ is obtained.
Thus we have a~solution $(Y_n, Z_n)$ $(n\in \mathbb {Z} )$ which satisf\/ies
the condition $(\frak{y}_{n_o}, Y_{n_o} )= (\tilde{\frak{y} }_o, y_o)$
and $(\frak{z}_{n_o}, Z_{n_o} )= ( \tilde{\frak{z} }_o, z_o)$.
\end{proof}

The uniqueness of the solution is not satisf\/ied for some cases.
In the case $U=U'$ or $V=V'$ in the proof of the proposition, i.e.\
the case $Y_m= A_3$, $A_4$, $A_1+mQ$ or $A_2+mQ$, we do not have uniqueness of the solution.
Similarly in the case $Z_{m+1}= B_3$, $B_4$, $B_1+mQ$ or $B_2+mQ$, we do not have uniqueness of the
solution.
On the other hand, in the case $\frak{y}_m =\frak{z}_m=-1$ for all~$m$, i.e.\
the case essentially equivalent to the ultradiscrete Painlev\'e VI without parity variables introduced by
Ormerod~\cite{Or} (see also~\cite{Tsu}), we have the uniqueness of the initial value problem (see
equations~\eqref{eq:om-zmzm1+}, \eqref{eq:om-ymym1+}).

\section{Ultradiscrete Riccati-type equation \\ and ultradiscrete Painlev\'e VI equation}
\label{sec:UDRUDP6}

We show that the ultradiscrete Riccati-type equation with parity variables
(equations~\eqref{eq:udRiccati2}, \eqref{eq:udRiccati1}) satisf\/ies the ultradiscrete Painlev\'e VI equation with
parity variables (equations~\eqref{eq:udP6zz}, \eqref{eq:udP6yy}).
\begin{lemma}
If $\max(X_1,X_2)=\max(X_3,X_4)$ and $\max(W_1,W_2)=\max(W_3,W_4)$, then
\begin{gather}
\max(X_1+W_1,X_3+W_3,X_2+W_4,X_4+W_2)
\nonumber
\\
\qquad
=\max(X_2+W_2,X_4+W_4,X_1+W_3,X_3+W_1).
\label{eq:lemeq}
\end{gather}
\end{lemma}

This lemma can be shown by case-by-case analysis (see~\cite{Tsu}, e.g.\
if $X_1 \geq X_2$, $X_3 \geq X_4$, $W_1 \leq W_2$ and $W_3 \leq W_4$, then $X_1=X_3$, $W_2=W_4$ and l.h.s.\
of equation~\eqref{eq:lemeq} $= \max(X_1+W_1,X_1+W_3,X_2+W_2,X_4+W_2) =$ r.h.s.\
of equation~\eqref{eq:lemeq}).
Note that it is an ultradiscrete analogue of the identity
$e^{X_1}-e^{X_3}=e^{X_4}-e^{X_2}$, $e^{W_1}-e^{W_3}=e^{W_4}-e^{W_2} \Longrightarrow
e^{X_1+W_1}+e^{X_3+W_3}+e^{X_2+W_4}+e^{X_4+W_2}=e^{X_2+W_2}+e^{X_4+W_4}+e^{X_2+W_3}+e^{X_3+W_1}$.
\begin{theorem}
If $A_1+B_3+Q=A_3+B_1$ and $A_2+B_4=A_4+B_2$, then any solutions of the ultradiscrete Riccati-type equation
with parity variables $($equations~\eqref{eq:udRiccati2}, \eqref{eq:udRiccati1}$)$ satisfy the ultradiscrete Painlev\'e
VI equation with parity variables $($equations~\eqref{eq:udP6zz}, \eqref{eq:udP6yy}$)$.
\end{theorem}
\begin{proof}
This theorem is shown by the lemma and case-by-case analysis.

We show equation~\eqref{eq:udP6yy} for every values of parity variables $\frak{z}_{m+1}$, $\frak{y}_{m}$, $\frak{y}_{m+1}$
by applying equations~\eqref{eq:udRiccati2}, \eqref{eq:udRiccati1} and the lemma.

In the case $\frak{z}_{m+1}= \frak{y}_{m} =\frak{y}_{m+1} =1$, we set $X_1= mQ+A_2+B_4$, $X_2=
Y_m+Z_{m+1}$, $X_3=A_4+Z_{m+1}$, $X_4=B_4+Y_m$, $W_1=mQ+A_3+B_1$, $W_2=Y_{m+1}+Z_{m+1}$, $W_3=A_3+Z_{m+1}$,
$W_4=B_3+Y_{m+1}$.
Then equations~\eqref{eq:udRiccati2}, \eqref{eq:udRiccati1} are written as $\max(X_1,X_2)=\max(X_3,X_4)$ and
$\max(W_1,W_2)=\max(W_3,W_4)$, i.e.\ equations~\eqref{eq:udR2++}, \eqref{eq:udR1++}.
By applying the lemma, we have
\begin{gather*}
\max\big[(mQ+A_2+B_4)+(mQ+A_3+B_1),(A_4+Z_{m+1})+(A_3+Z_{m+1}),
\\
\qquad
(Y_m+Z_{m+1})+(B_3+Y_{m+1}),(B_4+Y_m)+(Y_{m+1}+Z_{m+1})\big]
\\
\quad
=\max\big[(Y_m+Z_{m+1})+(Y_{m+1}+Z_{m+1}),(B_4+Y_m)+(B_3+Y_{m+1}),
\\
\qquad
(mQ+A_2+B_4)+(A_3+Z_{m+1}),(A_4+Z_{m+1})+(mQ+A_3+B_1)\big],
\end{gather*}
which is equivalent to equation~\eqref{eq:udP6yy} in the case $\frak{z}_{m+1}= \frak{y}_{m} =\frak{y}_{m+1}
=1$ (i.e.\ equation~\eqref{eq:zm1+omom1+}) by using the relation $A_2+B_4=A_4+B_2$.

We can show similarly equation~\eqref{eq:udP6yy} in the cases $(\frak{z}_{m+1}, \frak{y}_{m}, \frak{y}_{m+1})=
(1,1,-1)$, $(1,-1,1)$, $(1,-1,-1)$ by applying the lemma.

In the case $(\frak{z}_{m+1}, \frak{y}_{m}, \frak{y}_{m+1})= (-1,1,1)$, equation~\eqref{eq:udP6yy} is shown by
using equation~\eqref{eq:ass}.

In the cases $(\frak{z}_{m+1}, \frak{y}_{m}, \frak{y}_{m+1})= (-1,-1,1)$, $(-1,1,-1)$, $(-1,-1,-1)$, there
is no solution to equations~\eqref{eq:udRiccati2}, \eqref{eq:udRiccati1}, and the theorem holds true.

Similarly equation~\eqref{eq:udP6zz} is shown for every values of parity variables
$\frak{y}_{m+1}$, $\frak{z}_{m}$, $\frak{z}_{m+1}$ by applying equations~\eqref{eq:udRiccati2}, \eqref{eq:udRiccati1}
and using the relation
$A_1+B_3+Q=A_3+B_1 $.
For details see~\cite{Tsu}.
\end{proof}
The ultradiscrete Riccati equation in this paper never appear in the case $\frak{y}_m =\frak{z}_m=-1 $,
i.e.\
the case without parity variable discussed in~\cite{Or}.
Therefore ultradiscretization with parity variable is essential to obtain the ultradiscrete Riccati
equation.

\section{Solutions}
\label{sec:sol}
\subsection{Solutions of the ultradiscrete Riccati-type equation} We directly investigate solutions of the
ultradiscrete Riccati-type equation with parity variables (see
equations~\eqref{eq:udRiccati2}, \eqref{eq:udRiccati1}), which are also solutions of the ultradiscrete Painlev\'e
VI equation with parity variables (equations~\eqref{eq:udP6zz}, \eqref{eq:udP6yy}).
Although it would be natural to investigate solutions of the ultradiscrete Riccati-type equation by
introducing the ultradiscrete hypergeometric equation, we leave it to a~future problem.
In this subsection, we assume that $Q>0$, $A_1+B_3+Q=A_3+B_1$ and $A_2+B_4=A_4+B_2$.
Set
\begin{gather*}
h=A_3+B_1-A_2-B_4,
\qquad
h'=A_3-A_4-B_3+B_4.
\end{gather*}
To f\/ind solutions of the ultradiscrete Riccati-type equation with parity variables, we set an ansatz that
$\frak{y}_m=-1$ and $\frak{z}_{m}=+1$.
Then equations~\eqref{eq:udRiccati2}, \eqref{eq:udRiccati1} are written as
\begin{gather*}
\max(mQ+A_2,Y_m)+B_4=Z_{m+1}+\max(A_4,Y_m),
\\
\max((m+1)Q+A_1,Y_{m+1})+B_3=Z_{m+1}+\max(A_3,Y_{m+1}).
\end{gather*}
We specify the maximum in each term by
\begin{gather}
mQ+A_2\geq Y_m, \qquad A_4\leq Y_{m}, \qquad (m+1)Q+A_1\geq Y_{m+1}, \qquad A_3\leq Y_{m+1}.
\label{eq:ineqs1}
\end{gather}
Then we have $(m+1)Q+A_1 +B_3=Z_{m+1}+ Y_{m+1} $, $mQ+A_2 +B_4 =Z_{m+1}+Y_m $ and $Y_{m+1} -Y_m =A_3+B_1-
A_2-B_4 =h$.
Hence we obtain solutions $Y_m =hm+c $, $Z_{m+1} = (Q-h)m +A_2+B_4 -c$ where $c$ is a~constant which
satisf\/ies inequalities~\eqref{eq:ineqs1}.
Thus we have four solutions of equations~\eqref{eq:udRiccati2}, \eqref{eq:udRiccati1} which has a~parameter $c$ or
$c'$ with conditions as follows:
\begin{gather*}
(\frak{y}_m,Y_m)=(-1,hm+c), \qquad (\frak{z}_{m+1},Z_{m+1})=(+1,(Q-h)m+A_2+B_4-c),
\\
\text{condition for equation~\eqref{eq:udRiccati2}}:\qquad hm\geq A_4-c, \qquad (Q-h)m\geq c-A_2,
\\
\text{condition for equation~\eqref{eq:udRiccati1}}:\qquad h(m+1)\geq A_3-c, \qquad (Q-h)(m+1)\geq c-A_1,
\\
(\frak{y}_m,Y_m)=(+1,hm+A_2+B_4-c), \qquad (\frak{z}_{m+1},Z_{m+1})=(-1,(Q-h)m+c),
\\
\text{condition for equation~\eqref{eq:udRiccati2}}:\qquad (Q-h)m\geq B_4-c, \qquad hm\geq c-B_2,
\\
\text{condition for equation~\eqref{eq:udRiccati1}}:\qquad (Q-h)m\geq B_3-c, \qquad hm\geq c-B_1,
\\
(\frak{y}_m,Y_m)=(-1,h'm+c'), \qquad (\frak{z}_{m+1},Z_{m+1})=(+1,h'm+c'+B_4-A_4),
\\
\text{condition for equation~\eqref{eq:udRiccati2}}:\qquad h'm\leq A_4-c', \qquad (Q-h')m\leq c'-A_2,
\\
\text{condition for equation~\eqref{eq:udRiccati1}}:\qquad h'(m+1)\leq A_3-c', \qquad (Q-h')(m+1)\leq c'-A_1,
\\
(\frak{y}_m,Y_m)=(+1,h'm+c'-B_4+A_4), \qquad (\frak{z}_{m+1},Z_{m+1})=(-1,h'm+c'),
\\
\text{condition for equation~\eqref{eq:udRiccati2}}:\qquad h'm\leq B_4-c', \qquad (Q-h')m\leq c'-B_2,
\\
\text{condition for equation~\eqref{eq:udRiccati1}}:\qquad h'm\leq B_3-c', \qquad (Q-h')m\leq c'-B_1.
\end{gather*}
We also have solutions of equations~\eqref{eq:udRiccati2}, \eqref{eq:udRiccati1} which do not contain parameters.
One of them is the following solution with the conditions:
\begin{gather*}
(\frak{y}_{m+1},Y_{m+1})=(-1,A_2+B_4-B_1), \qquad (\frak{z}_{m+1},Z_{m+1})=(+1,mQ+B_1),
\\
A_2+B_4\leq A_3+B_1, \qquad B_1\leq B_2, \qquad mQ\geq\max(A_2+B_4-A_1-B_1,B_4-B_1).
\end{gather*}
Besides this type, we have four types of solutions for $m \gg 0$ in each case of the signs $(\frak{y}_{m},
\frak{z}_{m}) =(-1, +1)$, $(+1,-1)$ and $(+1,+1)$ respectively.
We also have four types of solutions for $m \ll 0$ in each case of the signs $(\frak{y}_{m}, \frak{z}_{m}) =(-1, +1)$, $(+1,-1)$ and $(+1,+1)$ respectively.
One of the solutions is as follows:
\begin{gather*}
(\frak{y}_{m+1},Y_{m+1})=(+1,A_3), \qquad (\frak{z}_{m+1},Z_{m+1})=(+1,B_4),
\\
A_3\geq A_4, \qquad B_3\leq B_4, \qquad mQ\leq\min(A_3-A_2,B_4-B_1).
\end{gather*}

We can construct global solutions of the ultradiscrete Riccati-type equation with parity variab\-les by
patching solutions for each region of the variable~$m$ suitably.
If there exists a~value~$c$ such that
\begin{gather*}
\max(A_1,A_4,A_2\!+\!B_4\!-\!B_1,A_1\!-\!B_1\!+\!B_2)\leq c\leq\min(A_2,A_3,A_3\!-\!B_3\!+\!B_4,A_2\!+\!B_4\!-\!B_3),
\end{gather*}
(the conditions $0\leq h \leq Q$ and $0\leq h' \leq Q$ are implied), we have a~solution written as
\begin{gather}
(\frak{y}_m,Y_m)=\begin{cases}
(-1,h'm+c),&m\leq0,
\\
(-1,hm+c),&m\geq1,
\end{cases}
\nonumber
\\
(\frak{z}_{m},Z_{m})=\begin{cases}
(+1,h'm+B_3-A_3+c),&m\leq0,
\\
(+1,(Q-h)m+A_1+B_3-c),&m\geq1.
\end{cases}
\label{eq:sol0}
\end{gather}
Let $m_0 <0$.
If $0\leq h\leq Q$, $0 \leq h' \leq Q$, $B_3\leq B_4 \leq B_1\leq B_4 +Q $, $A_3 +B_1 \geq A_4 +B_4$,
$\max( A_2,A_4 )\leq A_3$ and there exists a~value $c'$ such that
\begin{gather*}
A_4\leq h'm_0+c'\leq\min(A_3,A_4+h'), \qquad (Q-h')m_0+\max(A_1,A_2)\leq c',
\end{gather*}
then we have a~solution written as
\begin{gather}
(\frak{y}_m,Y_m)=\begin{cases}
(-1,h'm+c'),&m\leq m_0,
\\
(+1,A_3),&m_0+1\leq m\leq0,
\\
(-1,hm+A_2),&m\geq1,
\end{cases}
\nonumber
\\
(\frak{z}_{m},Z_{m})=\begin{cases}
(+1,h'm+B_3-A_3+c'),&m\leq m_0,
\\
(+1,B_4),&m_0+1\leq m\leq1,
\\
(+1,(Q-h)(m-1)+B_4),&m\geq2.
\end{cases}
\label{eq:soln2}
\end{gather}
Let $m_0 >0$.
If $0 \leq h \leq Q$, $0 \leq h' \leq Q$, $B_4 \leq B_1 \leq B_2$, $A_1 +B_1 \leq A_2 +B_4 \leq A_1+B_1
+Q$, $A_2 \leq \min (A_1,A_3)+Q$, $A_1 \geq A_4$, $A_2+B_3 \leq B_4+A_3+Q$ and there exists a~value $c$
such that
\begin{gather*}
h(m_0-1)+B_1\leq c\leq h m_0+\min(B_1,B_2), \qquad (Q-h)m_0+c\geq\max(B_3+(Q-h),B_4),
\end{gather*}
then we have a~solution written as
\begin{gather}
(\frak{y}_m,Y_m)=\begin{cases}
(+1,h'm+A_1-B_1+B_2),&m\leq-1,
\\
(-1,A_2+B_4-B_1),&0\leq m\leq m_0-1,
\\
(+1,hm+A_2+B_4-c),&m\geq m_0,
\end{cases}
\nonumber
\\
(\frak{z}_{m},Z_{m})=\begin{cases}
(-1,h'm+B_2-Q),&m\leq-1,
\\
(+1,A_2+B_4-A_1-Q),&m=0,
\\
(+1,(m-1)Q+B_1),&1\leq m\leq m_0,
\\
(-1,(Q-h)(m-1)+c),&m\geq m_0+1.
\end{cases}
\label{eq:solp}
\end{gather}
Note that we have other solutions which can be obtained similarly.

We now f\/ix the parameters by $A_1=25$, $A_2=46$, $A_3=67$, $A_4=23$, $B_1 = 59$, $B_2= 65$, $B_3=1$, $B_4=42$, $Q=100$.
Then $B_1 +A_3=A_1+B_3+Q$, $A_4+B_2=A_2+B_4$, $h= 38$ and $h'=85 $.
Equation~\eqref{eq:sol0} is a~solution of the ultradiscrete Riccati-type equation with parity variables, if $31
\leq c \leq 46$.
Moreover equation~\eqref{eq:soln2} (resp.\
equation~\eqref{eq:solp}) is a~solution of the ultradiscrete Riccati-type equation with parity variables, if $m_0
<0$ and $-85 m_0 +23 \leq c'\leq -85 m_0+ 67$ (resp.\
{}$m_0 >0$ and $ 21 + 38 m_0 \leq c\leq 59 + 38 m_0 $).
We also have other solutions.
For example, if $m_0 <0$ and $-85 m_0 -18\leq c'\leq -85 m_0+ 23$, then we have the following solution:
\begin{gather*}
(\frak{y}_m,Y_m)=\begin{cases}
(-1,85m+c'),&m\leq m_0,
\\
(+1,67),&m_0+1\leq m\leq0,
\\
(-1,38m+46),&m\geq1,
\end{cases}
\\
(\frak{z}_{m},Z_{m})=\begin{cases}
(+1,85m-66+c'),&m\leq m_0+1,
\\
(+1,42),&m_0+2\leq m\leq1,
\\
(+1,62m-20),&m\geq2.
\end{cases}
\end{gather*}
If $m_0 <0$ and $15 m_0 +31 \leq c'\leq 15m_0 +46 $, then we have the following solution:
\begin{gather*}
(\frak{y}_m,Y_m)=\begin{cases}
(-1,85m+c'),&m\leq m_0,
\\
(-1,100m+31),&m_0+1\leq m\leq0,
\\
(-1,38m+31),&m\geq1,
\end{cases}
\\
(\frak{z}_{m},Z_{m})=\begin{cases}
(+1,85m-66+c'),&m\leq m_0,
\\
(+1,100m-35),&m_0+1\leq m\leq0,
\\
(+1,62m-5),&m\geq1.
\end{cases}
\end{gather*}

\subsection{Solutions of the ultradiscrete Painlev\'e VI equation\\ without parity variables}

We now investigate solutions of ultradiscrete Painlev\'e VI with the f\/ixed parity variable $\frak{y}_m
=\frak{z}_m=-1$ for all $m$ (see equations~\eqref{eq:om-zmzm1+}, \eqref{eq:om-ymym1+}).

We look for the solutions written as $Y_m =\delta m +\beta $ and $Z_m = \alpha m +\gamma $ for $m \gg 0$.
We substitute them into equations~\eqref{eq:om-zmzm1+}, \eqref{eq:om-ymym1+}.
Then we have $\alpha + \delta =Q $, $2(\beta +\gamma)+ \alpha = B_3+B_4+A_1+A_2$ and inequalities among
the parameter.
More precisely, if
\begin{gather}
0\leq\alpha\leq Q, \qquad 2(\beta+\gamma)+\alpha=B_3+B_4+A_1+A_2,
\nonumber
\\
\alpha(m+1)+\gamma\geq\max(B_3,B_4), \qquad \alpha m+\min(A_1,A_2)\geq\beta,
\nonumber
\\
(Q-\alpha)m+\beta\geq\max(A_3,A_4), \qquad (Q-\alpha)m+\min(B_1,B_2)\geq\alpha+\gamma,
\label{eq:albega1}
\end{gather}
then the functions
\begin{gather*}
Y_m=(Q-\alpha)m+\beta,
\qquad
Z_m=\alpha m+\gamma
\end{gather*}
satisfy equations~\eqref{eq:om-zmzm1+}, \eqref{eq:om-ymym1+}.
Note that the constraint condition $ B_3+B_4+A_1+A_2=Q+ A_3+A_4+B_1+B_2$ is used to prove the equalities.
Similarly, if
\begin{gather}
0\leq\alpha'\leq Q, \qquad \alpha'+2(\gamma'-\beta')=B_3+B_4-A_3-A_4,
\nonumber\\
\alpha'(m+1)+\gamma'\leq\min(B_3,B_4), \qquad \alpha'm+\beta'\leq\min(A_3,A_4),
\nonumber\\
(Q-\alpha')m+\max(B_1,B_2)\leq\alpha'+\gamma', \qquad (Q-\alpha')m+\max(A_1,A_2)\leq\beta',
\label{eq:albega2}
\end{gather}
then the functions
\begin{gather*}
Y_m=\alpha'm+\beta',
\qquad
Z_m=\alpha'm+\gamma'
\end{gather*}
satisfy equations~\eqref{eq:om-zmzm1+}, \eqref{eq:om-ymym1+}.
We propose the following conjecture for solutions of ultradiscrete Painlev\'e VI without parity variables
(equations~\eqref{eq:om-zmzm1+}, \eqref{eq:om-ymym1+}), which is supported by several numerical solutions.
\begin{conj}
\label{conj:sol}
For every solution $Y_m$, $Z_m $ $(m \in \mathbb {Z})$ of ultradiscrete Painlev\'e VI without parity
variable $($equations~\eqref{eq:om-zmzm1+}, \eqref{eq:om-ymym1+}$)$, there exist $m_0, m_0' \in \mathbb {Z}$ and
$\alpha, \beta, \gamma, \alpha', \beta ', \gamma ' \in \mathbb {R}$ satisfying
equations~\eqref{eq:albega1}, \eqref{eq:albega2} such that
\begin{gather*}
Y_m=\alpha'm+\beta', \qquad Z_m=\alpha'm+\gamma', \qquad m\leq m'_0,
\\
Y_m=(Q-\alpha)m+\beta, \qquad Z_m=\alpha m+\gamma, \qquad m\geq m_0.
\end{gather*}
\end{conj}
We give examples of solutions.
We consider solutions in the case $A_1=32$, $A_2=33$, $A_3=37$, $A_4=22$, $B_1=53$, $B_2=65$, $B_3=8$, $B_4=4$, and $
Q=100$.
We choose the initial values by $Y_0=43$, $Z_0=40$.
Then the solution is written as
\begin{gather*}
Y_m=\begin{cases}
95m+111,&m\leq-1,
\\
43,&m=0,
\\
11m+111,&m\geq1,
\end{cases}
\qquad
Z_m=\begin{cases}
95m+40,&m\leq0,
\\
89m-117,&m\geq1.
\end{cases}
\end{gather*}
On dif\/ferent initial values $Y_0=43$, $Z_0=50$, the solution is
\begin{gather*}
Y_m=\begin{cases}
85m-81,&m\leq-8,
\\
-669,&m=-7,
\\
115m+131,&-1\leq m\leq-6,
\\
43,&m=0,
\\
-9m+131,&1\leq m\leq11,
\\
9m-72,&m\geq12,
\end{cases}
\qquad
Z_m=\begin{cases}
85m-147,&m\leq-7,
\\
115m+50,&0\leq m\leq-6,
\\
109m-147,&1\leq m\leq11,
\\
1156,&m=12,
\\
91m+65,&m\geq13.
\end{cases}
\end{gather*}
In these cases, the conjecture is true.

\section{Concluding remarks}

In Section~\ref{sec:UD}, we obtained a~ultradiscretization with parity
variables of the $q$-dif\/ference Pain\-le\-v\'e~VI equation.
A list of Painlev\'e-type equations of second order was obtained by Sakai~\cite{Sak0}, and some members in
the list are $q$-dif\/ference Painlev\'e equations.
We believe that ultradiscretization with parity variables of the $q$-dif\/ference Painlev\'e equations can be done.

Although we investigated solutions of ultradiscrete Riccati-type equation directly, we did not study
ultradiscrete hypergeometric equation in this paper.
A theory of ultradiscrete hypergeometric equations should be developed because it will have potential for
applications to several equations including $q$-dif\/ference hypergeometric equations and ultradiscrete
Painlev\'e equations.

A merit of ultradiscrete equations is that we may have exact solutions with the aid of com\-pu\-ters.
We formulated Conjecture~\ref{conj:sol} by calculating several solutions in use of a~computer.
On the other hand, Murata~\cite{Mu} obtained exact solutions with two parameters for a~ultradiscrete
Painlev\'e~II equation.
We hope to understand exact solutions for ultradiscrete Painlev\'e equations deeply.

\subsection*{Acknowledgments} The authors would like to thank Professor Junkichi Satsuma for discussions
and suggestions.
They also thank the referees for valuable comments.
The f\/irst author is partially supported by the Grant-in-Aid for Young Scientists~(B) (No.~22740107) from
the Japan Society for the Promotion of Science.

\pdfbookmark[1]{References}{ref}
\LastPageEnding

\end{document}